\documentclass[11pt]{article}
\usepackage[dvipdfmx]{graphicx}
\usepackage{amsbsy,amsfonts,amsmath,amssymb,amsthm,ascmac,mathrsfs,centernot}
\usepackage{graphics,epsfig,enumerate,fancybox,centernot,eclbkbox,color,framed}
\definecolor{shadecolor}{gray}{0.85}
\usepackage{url}
\usepackage{hyperref}
\hypersetup{colorlinks,citecolor=blue,linkcolor=blue,urlcolor=black}
\usepackage{bbm}
\newcommand{\indic}{\mathbbm{1}}
\usepackage{multirow}
\usepackage{rotating}
\usepackage{subcaption}
\definecolor{mygreen}{RGB}{0,128,0}
\DeclareMathOperator*{\argmin}{arg\,min}

\allowdisplaybreaks[1]

\theoremstyle{plain}
\numberwithin{equation}{section}
\newtheorem{thm}{Theorem}[section]

\newtheorem{prp}[thm]{Proposition}
\newtheorem{rmk}[thm]{Remark}

\newcommand{\dpst}{\displaystyle}
\newcommand{\EG}{E^{\sss\rm G}}
\newcommand{\ESCA}{E^{\sss\rm SCA}}

\newcommand{\ind}[1]{\indic_{\{#1\}}}
\newcommand{\lbeq}[1]{\label{eq:#1}}
\newcommand{\mE}{{\mathbb E}}
\newcommand{\N}{{\mathbb N}}
\newcommand{\nn}{\nonumber}

\newcommand{\mP}{{\mathbb P}}

\newcommand{\PG}{P^{\sss\rm G}}
\newcommand{\PSCA}{P^{\sss\rm SCA}}
\newcommand{\piG}{\pi^{\sss\rm G}}
\newcommand{\piSCA}{\pi^{\sss\rm SCA}}

\newcommand{\refeq}[1]{(\ref{eq:#1})}
\newcommand{\rhoTM}{\rho_{\sss\rm TM}}
\newcommand{\sss}{\scriptscriptstyle}

\newcommand{\tSCA}{t_{\rm mix}^{\sss\rm SCA}}

\newcommand{\vep}{\varepsilon}

\newcommand{\wG}{w^{\sss\rm G}}
\newcommand{\wSCA}{w^{\sss\rm SCA}}

\newcommand{\Evec}  {\textbf{E}}

\newcommand{\Pvec}  {\textbf{P}}

\newcommand{\qvec}  {\boldsymbol{q}}

\newcommand{\etavec}  {\boldsymbol{\eta}}
\newcommand{\sigmavec}  {\boldsymbol{\sigma}}
\newcommand{\tauvec}  {\boldsymbol{\tau}}
\newcommand{\xivec}  {\boldsymbol{\xi}}

\title{Mixing time and simulated annealing for\\
the stochastic cellular automata}
\author{
Bruno Hideki Fukushima-Kimura\footnote{Faculty of Science, Hokkaido University.}\hskip5pt${}^,$\footnote{Corresponding author: bruno@math.sci.hokudai.ac.jp}\\
Satoshi Handa\footnote{atama plus K.K.}\\
Katsuhiro Kamakura\footnote{TOME R\&D Inc.}\\
Yoshinori Kamijima\footnote{National Center for Theoretical Sciences (NCTS).}\\
Kazushi Kawamura\footnote{Institute of Innovative Research, Tokyo Institute of Technology.}\\
Akira Sakai\footnotemark[1]\hskip5pt${}^,$\footnote{\url{https://orcid.org/0000-0003-0943-7842}}
}

\begin{document}
\maketitle

\begin{abstract}
Finding a ground state of a given  Hamiltonian of an Ising model on a graph $G=(V,E)$ is an 
important but hard problem.  The standard approach for this kind of problem is the application of algorithms that rely on 
single-spin-flip Markov chain Monte Carlo methods, such as the simulated annealing based on Glauber or Metropolis dynamics.  In this paper, 
we investigate a particular kind of stochastic cellular automata, in which all spins are 
updated independently and simultaneously.  We prove that (i) if the 
temperature is fixed sufficiently high, then the mixing time is at 
most of order $\log|V|$, and that (ii) if the temperature drops in time 
$n$ as $1/\log n$, then the limiting measure is uniformly distributed 
over the ground states. We also provide some simulations of the algorithms studied in this paper implemented on a GPU and show their superior performance compared to the conventional simulated annealing.
\end{abstract}

\section{Introduction and main results}
There are several occasions in real life when we have to choose one
among extremely many options quickly.  In addition, we want our choice to be optimal in 
a certain sense. The so-called combinatorial optimization problems \cite{Lawler1976,Papa82} are ubiquitous and 
possibly quite hard to be solved in a fast way.  {In particular, polynomial-time algorithms for NP-hard problems are not known to exist~\cite{gj79}.}

A possible approach as an attempt to provide an optimal solution for a given problem would be to translate it
into the problem of minimizing the Hamiltonian of an Ising model on a finite graph such that each of its ground states corresponds to an optimal solution to the original problem and vice versa. See, for instance, \cite{l14} for a broad list of examples of such mappings. Let us consider a finite graph $G=(V,E)$ with no multi- or self-edges.
Given a collection of spin-spin coupling constants $\{J_{x,y}\}_{x,y \in V}$ which is symmetric (that is, $J_{x,y} = J_{y,x}$ for all $x,y$) and satisfies
$J_{x,y}=0$ whenever $\{x,y\}\notin E$, and a collection of external fields 
$\{h_x\}_{x\in V}$, the Hamiltonian of an Ising spin configuration 
$\sigmavec=\{\sigma_x\}_{x\in V}\in\Omega\equiv\{\pm1\}^V$ is defined by
\begin{align}
H(\sigmavec)=-\sum_{\{x,y\}\in E}J_{x,y}\sigma_x\sigma_y-\sum_{x\in V}h_x
 \sigma_x\equiv-\frac12\sum_{x,y\in V}J_{x,y}\sigma_x\sigma_y-\sum_{x\in V}h_x
 \sigma_x.
\end{align}
Let GS denote the set of its ground states, the configurations at which the Hamiltonian attains its 
minimum value, i.e., 
\begin{align}
\mathrm{GS}=\underset{\sigmavec}{\arg\min}H(\sigmavec)\equiv\big\{\sigmavec\in\Omega:
 H(\sigmavec)=\min_{\tauvec}H(\tauvec)\big\}.
\end{align}

A method that can possibly be considered in the search for ground states consists of using a Markov chain Monte Carlo 
(MCMC) to sample the Gibbs distribution $\piG_\beta\propto e^{-\beta H}$ at 
the inverse temperature $\beta\ge0$ given by
\begin{align}\lbeq{piGdef}
\piG_\beta(\sigmavec)=\frac{\wG_\beta(\sigmavec)}{\sum_{\tauvec}
 \wG_\beta(\tauvec)},\qquad\text{where}\quad
\wG_\beta(\sigmavec)=e^{-\beta H(\sigmavec)}.
\end{align}
It is straightforward to verify that the Gibbs distribution reaches its highest peaks on GS.
There are several MCMCs that can be applied to generate the Gibbs distribution as the 
equilibrium measure, one of them is the Glauber dynamics~\cite{g63}, 
whose transition matrix $\PG_\beta$ is defined by  
\begin{align}\lbeq{PGdef}
\PG_\beta(\sigmavec,\tauvec)=
 \begin{cases}
 \dpst\frac1{|V|}\frac{\wG_\beta(\sigmavec^x)}{\wG_\beta(\sigmavec)+\wG_\beta
  (\sigmavec^x)}&[\tauvec=\sigmavec^x],\\[1pc]
 \dpst1-\sum_{x\in V}\PG_\beta(\sigmavec,\sigmavec^x)&[\tauvec=\sigmavec],\\
 0&[\text{otherwise}],
 \end{cases}\qquad\text{where}\quad
(\sigmavec^x)_y=
 \begin{cases}
 \sigma_y&[y\ne x],\\
 -\sigma_y&[y=x].
 \end{cases}
\end{align}
Notice that, by introducing the cavity fields
\begin{align}\lbeq{cavity}
\tilde h_x(\sigmavec)=\sum_{y\in V}J_{x,y}\sigma_y+h_x,
\end{align}
the transition probability $\PG_\beta(\sigmavec,\sigmavec^x)$ can also be 
written as
\begin{align}\lbeq{PGrewr}
\PG_\beta(\sigmavec,\sigmavec^x)=\frac1{|V|}\frac{e^{-\beta\tilde h_x(\sigmavec)
 \sigma_x}}{2\cosh(\beta\tilde h_x(\sigmavec))}.
\end{align}
The transition probability above can be interpreted as the probability of choosing 
the vertex $x$ uniformly at random from $V$ and then flipping its spin value with probability 
proportional to $\wG_\beta(\sigmavec^x)$.  
Furthermore, since $\PG_\beta$ is aperiodic, irreducible, and reversible with respect to 
$\piG_\beta$ (i.e., the identity $\piG_\beta(\sigmavec)\PG_\beta(\sigmavec,\tauvec)
=\piG_\beta(\tauvec)\PG_\beta(\tauvec,\sigmavec)$ holds for all 
$\sigmavec,\tauvec\in\Omega$), it follows that the Glauber dynamics converges to its unique equilibrium 
distribution, namely, the Gibbs distribution $\piG_\beta$.  

In practice, the method that actually has been widely employed in several real-world applications is known as the simulated annealing algorithm \cite{AK1989,c85,IZRT15,kgv83,LA1987,Wong1988}, which motivated the main subject of this paper. The standard approach consists of applying a discrete-time inhomogeneous Markov chain based on single-spin-flip dynamics (such as Glauber dynamics or Metropolis dynamics \cite{Metropolis53}) where the temperature drops every time a spin value is updated. If the temperature is set to decrease at an appropriate rate, then it is guaranteed that such a procedure will asymptotically converge to one of the ground states, see \cite{AK1989,b99,h88}.

Note that, in the methods described above, the number of spin-flips 
per update is at most one, so, in principle, we may take some benefit if we consider methods that allow a larger number of spin-flips. In respect of ferromagnetic spin systems, the Swendsen-Wang algorithm 
\cite{sw87} is a cluster-flip Markov chain, in which many spins can be flipped 
simultaneously, unlike in the Glauber and other single-spin-flip dynamics.  
However, forming a cluster to be flipped yields strong dependency among 
spin variables.  

{In recent studies such as in \cite{osky19,STATICA}, some algorithms 
that rely on parallel and independent spin-flips have shown 
significantly higher performance in approximating ground states compared to some of the well-established algorithms based on single-spin-flip dynamics. For that reason,
such algorithms of that nature deserve some attention and a rigorous treatment  from the mathematical point of
view, in order to understand their mechanisms and limitations, becomes necessary. Furthermore, due to their main feature, the possibility of employing hardware accelerators such as annealing processors strongly relying on parallelization \cite{DA2019,IEEE2023,STATICA} to speed up simulations has a great appeal and may be advantageous for the computation of solutions of real-time and time-constrained problems. So, in Section \ref{sec:eSCA} we include some analyses  comparing the accuracy and the simulation times of Glauber dynamics and the algorithms present in this paper when implemented on a CPU and a GPU.}

In this paper, we investigate a particular class of probabilistic cellular automata, or PCA, 
studied in \cite{dss12,st18}.  Since the acronym PCA has already been long used 
to stand for principal component analysis in statistics, we would 
rather use the term stochastic cellular automata (SCA). Let us start by considering the so-called pinning parameters 
$\qvec=\{q_x\}_{x\in V}$ and using them to introduce 
an extended version of the Hamiltonian, whose expression is given by
\begin{align}\lbeq{extendedH}
\tilde H(\sigmavec,\tauvec)&=-\frac12\sum_{x,y\in V}J_{x,y}\sigma_x\tau_y
 -\frac12\sum_{x \in V}h_x(\sigma_x+\tau_x) - \frac{1}{2} \sum_{x \in V} q_{x}
 \sigma_{x} \tau_{x}\nn\\
&= -\frac12\sum_{x\in V}\big(\tilde h_x(\sigmavec)+q_x\sigma_x\big)\tau_x
 -\frac12\sum_{x\in V}h_x\sigma_x.
\end{align}
Note that the relationship between $H$ and $\tilde H$ is given by
\begin{align}\lbeq{symmetry}
\tilde H(\sigmavec,\sigmavec)=H(\sigmavec)-\frac12\sum_{x\in V}q_x.
\end{align}
Then, we introduce
\begin{align}
\wSCA_{\beta,\qvec}(\sigmavec)=\sum_{\tauvec}e^{-\beta \tilde H(\sigmavec,
 \tauvec)}\stackrel{\text{\refeq{extendedH}}}=\prod_{x\in V}2e^{\frac\beta2
 h_x\sigma_x}\cosh\big(\tfrac\beta2(\tilde h_x(\sigmavec)+q_x\sigma_x)\big),
\end{align}
and define the SCA transition matrix $\PSCA_{\beta,\qvec}$ by letting
\begin{align}\lbeq{PSCAdef}
\PSCA_{\beta,\qvec}(\sigmavec,\tauvec)=\frac{e^{-\beta \tilde H(\sigmavec,
 \tauvec)}}{\wSCA_{\beta,\qvec}(\sigmavec)}\stackrel{\text{\refeq{extendedH}}}
 = \prod_{x\in V}\frac{e^{\frac\beta2(\tilde h_x(\sigmavec)+q_x\sigma_x)
 \tau_x}}{2\cosh(\frac\beta2(\tilde h_x(\sigmavec)+q_x\sigma_x))}.
\end{align}
Due to the rightmost expression from equation \refeq{PSCAdef}, we can interpret that all spins in the system 
are updated independently and simultaneously according to certain local probability rules.  This implies that the SCA is allowed to 
move from any spin configuration to another in just one step, which, in principle, may potentially 
result in a faster convergence to equilibrium.  Since $\tilde H$ is 
symmetric (due to the symmetry of the spin-spin coupling constants), i.e.,
$\tilde H(\sigmavec,\tauvec)=\tilde H(\tauvec,\sigmavec)$, then,
the middle expression in equation \refeq{PSCAdef} implies that $\PSCA_{\beta,\qvec}$ is 
reversible with respect to the equilibrium distribution $\piSCA_{\beta,\qvec}$ given by
\begin{align}\lbeq{piSCAdef}
\piSCA_{\beta,\qvec}(\sigmavec)&=\frac{\wSCA_{\beta,\qvec}(\sigmavec)}
 {\sum_{\tauvec}\wSCA_{\beta,\qvec}(\tauvec)}.
\end{align}
Although this does not necessarily coincide with the Gibbs distribution, and therefore we cannot naively 
use it to search for the ground states, the total-variation distance (cf., 
\cite[Definition~4.1.1 \& (4.1.5)]{b99})
\begin{align}\lbeq{TVdistance}
\|\piSCA_{\beta,\qvec}-\piG_\beta\|_{\sss\rm TV}\equiv\frac12\sum_{\sigmavec}
 |\piSCA_{\beta,\qvec}(\sigmavec)-\piG_\beta(\sigmavec)|=1-\sum_{\sigmavec}
 \piSCA_{\beta,\qvec}(\sigmavec)\wedge\piG_\beta(\sigmavec)
\end{align}
tends to zero as $\min_xq_x\uparrow\infty$.  This is a positive aspect of 
the SCA with large $\qvec$.  On the other hand, since the off-diagonal entries 
of the transition matrix $\PSCA_{\beta,\qvec}$ tends to zero as 
$\min_xq_x\uparrow\infty$, the SCA with large $\qvec$ may well be much slower 
than expected.  This is why we call $\qvec$ the pinning parameters.

Having in mind the development of a simulated annealing algorithm based on  SCA aiming at determining the ground states of $H$,  in Section \ref{sec:SCAannealing}
we investigate the SCA with the set of pinning parameters $\qvec$ satisfying 
\begin{align}\lbeq{MA-pinning-parameter}
q_x\ge
 \begin{cases}
 \dpst\sum_{y\in V}|J_{x, y}|-\frac12\sum_{y\in C}|J_{x, y}|\quad
  &[x\in C],\\[1.5pc]
 \dpst\frac\lambda2 &[x\notin C],
 \end{cases}
\end{align}
where $C$ is an arbitrary subset of $V$ and $\lambda$ is the largest 
eigenvalue of the matrix $[-J_{x, y}]_{V\times V}$.  This is a sufficient 
condition on $\tilde H$ that assures that its minimum value is attained on the diagonal entries, that is,
\begin{equation}
\min_{\sigmavec,\tauvec\in\Omega}\tilde H(\sigmavec,\tauvec)
 =\min_{\sigmavec\in\Omega}\tilde H(\sigmavec,\sigmavec).
\end{equation}
Condition \refeq{MA-pinning-parameter} originated in \cite{osky19}, where on its supplemental material a rather more general result is proved, establishing that, under this assumption, the inequality  
\begin{equation}\lbeq{mindiag}
    \tilde H(\sigmavec,\tauvec) > \min_{\sigmavec' \in \Omega} \tilde H(\sigmavec',\sigmavec')
\end{equation}
holds whenever $\sigmavec$ and $\tauvec$ are distinct, 
which, according to equation \refeq{symmetry}, implies that
\begin{equation}
    \argmin_{\sigmavec, \tauvec \in \Omega} \tilde H(\sigmavec,\tauvec) = \{(\sigmavec,\sigmavec): \sigmavec \in \mathrm{GS}\}.
\end{equation}
In this paper, we prove the following two statements:
\begin{enumerate}
\item[(i)]
If $\beta$ is sufficiently small and fixed, then the time-homogeneous SCA has 
a mixing time at most of order $\log|V|$ (Theorem~\ref{thm:mixing}).
\item[(ii)]
If $\beta_n$ increases in time $n$ as $\propto\log n$, then the 
time-inhomogeneous SCA weakly converges to $\piG_\infty$, the uniform distribution 
 over GS (Theorem~\ref{thm:SAforSCA}).
\end{enumerate}
The former result implies faster mixing than conventional single-spin-flip MCs, 
such as the Glauber dynamics (see Remark~\ref{rmk:mixing}(i)).  The latter 
refers to the applicability of the standard temperature-cooling schedule in the 
simulated annealing (see Remark~\ref{rmk:SAforSCA}(i)).

Although the two results above are proven mathematically rigorously, they may be difficult to be directly applied in practice.  
As mentioned earlier, the SCA is allowed to flip multiple spins in a single update, therefore,
in principle,  it potentially can reduce the mixing time compared to other single-spin-flip algorithms.  However, to 
attain such a small mixing time as in (i), we have to keep the temperature 
very high (as comparable to the radius of convergence of the high-temperature 
expansion, see \refeq{mixingcond} below).  Also, if we want to find a ground 
state by using an SCA-based simulated annealing algorithm, as stated in (ii), the 
temperature has to drop so slowly as $1/\log n$ (with a large multiplicative 
constant $\Gamma$, see \refeq{SCA-cooling-schedule} below), and therefore 
the number of steps required to come close to a ground state may well be 
extremely large.  The problem seems to be due to the introduction of the total-variation 
distance.  In order to make the distance $\|\mu-\nu\|_{\sss\rm TV}$ small, the two 
measures $\mu$ and $\nu$ have to be very close at every spin configuration.  
Moreover, since $\|\piSCA_{\beta,\qvec}-\piG_\beta\|_{\sss\rm TV}$ is not 
necessarily small for finite $\qvec$, we cannot tell anything about the excited 
states $\Omega\setminus\mathrm{GS}$.  In other words, we might be able to 
use the SCA under the condition \refeq{MA-pinning-parameter} to find the best 
options in combinatorial optimization, but not the second- or third-best options.  
To overcome those difficulties, we will shortly discuss a potential replacement 
for the total-variation distance at the end of Section \ref{sec:SCAannealing}.

In practice, to avoid applying logarithmic cooling schedules and having to wait for long execution times, faster cooling is often considered. In Section \ref{sec:eSCA}, we introduce a variant of the SCA called $\vep$-SCA and show in a series of examples that simulated annealing with exponential cooling schedules can be successfully applied to the Glauber dynamics, SCA, and $\vep$-SCA. The results are consistent with the preliminary findings from \cite{SSS21,TISCIE23}, which reveal that the SCA outperforms Glauber dynamics in most of the scenarios, while $\vep$-SCA outperforms both algorithms in all tested models. Providing a mathematical foundation for the $\vep$-SCA and proving a generalization of (ii) for exponential schedules comprise a future direction of our research.

\section{Mixing time for the SCA}\label{sec:mixing-time}
In this section, we show that the mixing in the SCA is faster than in the 
Glauber dynamics when the temperature is sufficiently high.  To do so, 
we first introduce some notions and notation. 

For $\sigmavec,\tauvec\in\Omega$, we let $D_{\sigmavec,\tauvec}$ be 
the set of vertices at which $\sigmavec$ and $\tauvec$ disagree:
\begin{align}
D_{\sigmavec,\tauvec}=\{x\in V:\sigma_x\ne\tau_x\}.
\end{align}
For a time-homogeneous Markov chain, whose $t$-step distribution $P^t$ 
converges to its equilibrium $\pi$, we define the mixing time 
as follows: given an $\vep\in[0,1]$, 
\begin{equation*}
t_{\rm mix}(\varepsilon)=\inf\Big\{t\ge0:\max_{\sigmavec}\|P^t(\sigmavec,\cdot)
 -\pi\|_{\sss\rm TV}\le\vep\Big\}.
\end{equation*}
In particular, we denote it by 
$\tSCA(\vep)$ when $P=\PSCA_{\beta,\qvec}$ and $\pi=\piSCA_{\beta,\qvec}$.
Then we define the transportation metric $\rhoTM$ between two probability measures on 
$\Omega$ as
\begin{equation}\lbeq{transporation-metric}
\rhoTM(\mu,\nu)=\inf\Big\{\Evec_{\mu,\nu}[|D_{X,Y}|]:
 \text{$(X, Y)$ is a coupling of $\mu$ and $\nu$}\Big\},
\end{equation}
where $\Evec_{\mu,\nu}$ is the expectation against the coupling measure 
$\Pvec_{\mu,\nu}$ whose marginals are $\mu$ for $X$ and $\nu$ for $Y$, 
respectively.
By \cite[Lemma~14.3]{lp17}, $\rhoTM$ indeed satisfies the axioms of metrics, 
in particular the triangle inequality: 
$\rhoTM(\mu,\nu)\le\rhoTM(\mu,\lambda)+\rhoTM(\lambda,\nu)$ 
holds for all probability measures $\mu,\nu,\lambda$ on $\Omega$.

The following is a summary of \cite[Theorem~14.6 \& Corollary~14.7]{lp17}, 
but stated in our context.

\begin{shaded}
\begin{prp}\label{prp:coupling}
If there is an $r\in(0,1)$ such that 
$\rhoTM(P(\sigmavec,\cdot),P(\tauvec,\cdot))\le r|D_{\sigmavec, \tauvec}|$ for 
any $\sigmavec, \tauvec \in \Omega$, then
\begin{align}\lbeq{upper-bound-dSCA}
\max_{\sigmavec\in\Omega}\|P^t(\sigmavec,\cdot)-\pi\|_{\sss\rm TV}\le r^t
 \max_{\sigmavec, \tauvec\in\Omega} |D_{\sigmavec, \tauvec}|.
\end{align}
Consequently,
\begin{align}\lbeq{upper-bound-SCA-mixing}
t_{\rm mix}(\vep)\le\bigg\lceil\frac{\log|V|-\log\vep}{\log(1/r)}\bigg\rceil.
\end{align}
\end{prp}
\end{shaded}

It is crucial to find a coupling $(X,Y)$ in which the size of $D_{X,Y}$ is 
decreasing in average, as stated in the hypothesis of the above proposition.  
Here is the main statement on the mixing time for the SCA.

\begin{shaded}
\begin{thm}\label{thm:mixing}
For any non-negative $\qvec$, if $\beta$ is sufficiently small such that, 
independently of $\{h_x\}_{x\in V}$,
\begin{align}\lbeq{mixingcond}
r\equiv\max_{x\in V}\bigg(\tanh\frac{\beta q_x}2+\sum_{y\in V}\tanh\frac{\beta
 |J_{x,y}|}2\bigg)<1,
\end{align}
then 
$\rhoTM(\PSCA_{\beta,\qvec}(\sigmavec,\cdot),\PSCA_{\beta,\qvec}(\tauvec,\cdot))
\le r|D_{\sigmavec, \tauvec}|$ for all $\sigmavec, \tauvec \in \Omega$, 
and therefore $\tSCA(\vep)$ obeys \refeq{upper-bound-SCA-mixing}.
\end{thm}
\end{shaded}

\begin{proof}
It suffices to show $\rhoTM(\PSCA_{\beta,\qvec}(\sigmavec,\cdot),\PSCA_{\beta,
\qvec}(\tauvec,\cdot))\le r$ for all $\sigmavec,\tauvec\in\Omega$ with 
$|D_{\sigmavec,\tauvec}|=1$.  If $|D_{\sigmavec,\tauvec}|\ge2$, then, by 
the triangle inequality along any sequence 
$(\etavec_0,\etavec_1,\dots,\etavec_{|D_{\sigmavec,\tauvec}|})$ of spin 
configurarions that satisfy
$\etavec_0=\sigmavec$, $\etavec_{|D_{\sigmavec,\tauvec}|}=\tauvec$ and 
$|D_{\etavec_{j-1},\etavec_j}|=1$ for all $j=1,\dots,|D_{\sigmavec,\tauvec}|$, 
we have
\begin{align}
\rhoTM\Big(\PSCA_{\beta,\qvec}(\sigmavec,\cdot),\PSCA_{\beta,\qvec}(\tauvec,
 \cdot)\Big)\le\sum_{j=1}^{|D_{\sigmavec,\tauvec}|}\rhoTM\Big(\PSCA_{\beta,
 \qvec}(\etavec_{j-1},\cdot),\PSCA_{\beta,\qvec}(\etavec_j,\cdot)\Big)\le r
 |D_{\sigmavec,\tauvec}|.
\end{align}

Suppose that $D_{\sigmavec,\tauvec}=\{x\}$, i.e., $\tauvec=\sigmavec^x$.  
For any $\sigmavec\in\Omega$ and 
$y\in V$, we let $p(\sigmavec,y)$ be the conditional SCA probability of 
$\sigma_y\to1$ given that the others are fixed (cf., \refeq{PSCAdef}):
\begin{align}\lbeq{condprob}
p(\sigmavec,y)=\frac{e^{\frac\beta2(\tilde h_y(\sigmavec)+q_y\sigma_y)}}
 {2\cosh(\frac\beta2(\tilde h_y(\sigmavec)+q_y\sigma_y))}=\frac{1+\tanh(
 \frac\beta2(\tilde h_y(\sigmavec)+q_y\sigma_y))}2.
\end{align}
Notice that $p(\sigmavec,y)\ne p(\sigmavec^x,y)$ only when $y=x$ or 
$y\in N_x\equiv\{v\in V:J_{x,v}\ne0\}$.  Using this as a threshold function for 
i.i.d.~uniform random variables $\{U_y\}_{y\in V}$ on $[0,1]$, we define the 
coupling $(X,Y)$ of $\PSCA_{\beta,\qvec}(\sigmavec,\cdot)$ and 
$\PSCA_{\beta,\qvec}(\sigmavec^x,\cdot)$ as
\begin{align}\lbeq{coupling-SCA}
X_y=\begin{cases}
 +1 & [U_y\le p(\sigmavec,y)],\\
 -1 & [U_y>p(\sigmavec,y)],
\end{cases}&&
Y_y=\begin{cases}
 +1 & [U_y\le p(\sigmavec^x,y)],\\
 -1 & [U_y>p(\sigmavec^x,y)].
\end{cases}
\end{align}
Denote the measure of this coupling by $\Pvec_{\sigmavec,\sigmavec^x}$ and 
its expectation by $\Evec_{\sigmavec,\sigmavec^x}$.  Then we obtain
\begin{gather}
\Evec_{\sigmavec,\sigmavec^x}[|D_{X,Y}|]=\Evec_{\sigmavec,\sigmavec^x}\bigg[
 \sum_{y\in V}\ind{X_y\ne Y_y}\bigg]=\sum_{y\in V}\Pvec_{\sigmavec,\sigmavec^x}
 (X_y\ne Y_y)=\sum_{y\in V}|p(\sigmavec,y)-p(\sigmavec^x,y)|\nn\\
=|p(\sigmavec,x)-p(\sigmavec^x,x)|+\sum_{y\in N_x}|p(\sigmavec,y)-p(\sigmavec^x,
 y)|,\lbeq{expected-hamming-SCA}
\end{gather}
where, by using the rightmost expression in \refeq{condprob},
\begin{align}
|p(\sigmavec,x)-p(\sigmavec^x,x)|\le\frac12\bigg|\tanh\bigg(\frac{\beta\tilde
 h_x(\sigmavec)}2+\frac{\beta q_x}2\bigg)-\tanh\bigg(\frac{\beta\tilde h_x
 (\sigmavec)}2-\frac{\beta q_x}2\bigg)\bigg|,
\end{align}
and for $y\in N_x$, 
\begin{align}
|p(\sigmavec,y)-p(\sigmavec^x,y)|\le\frac12\bigg|\tanh\bigg(\frac{\beta(\sum_{v
 \ne x}J_{v,y}\sigma_v+h_y+q_y\sigma_y)}2&+\frac{\beta J_{x,y}}2\bigg)\nn\\
-\tanh\bigg(\frac{\beta(\sum_{v\ne x}J_{v,y}\sigma_v+h_y+q_y\sigma_y)}2
 &-\frac{\beta J_{x,y}}2\bigg)\bigg|.
\end{align}
Since $|\tanh(a+b)-\tanh(a-b)|\le2\tanh|b|$ for any $a,b$, we can conclude 
\begin{align}
\rhoTM\Big(\PSCA_{\beta,\qvec}(\sigmavec,\cdot),\PSCA_{\beta,\qvec}(\sigmavec^x,
 \cdot)\Big)\le\Evec_{\sigmavec,\sigmavec^x}[|D_{X,Y}|]\le\tanh\frac{\beta q_x}2
 +\sum_{y\in N_x}\tanh\frac{\beta|J_{x,y}|}2\le r,
\end{align}
as required.
\end{proof}

\begin{rmk}\label{rmk:mixing}
{\rm
\begin{enumerate}[(i)]
\item
{It is known that, in a very general setting, the mixing time for the Glauber dynamics \refeq{PGdef}
is at least of order $|V| \log |V|$, see \cite{HS07}.  Therefore, Theorem~\ref{thm:mixing} implies that the 
SCA reaches equilibrium way faster than the Glauber dynamics, as long as the temperature is high enough. Even though the result above
does not play a role in practical applications and the SCA cannot be used to sample the Gibbs distribution, it does give us a hint that we may extract some benefit by considering multiple spin-flip algorithms to speed up simulations.}
\item
It is of some interest to investigate the average number of spin-flips per 
update, although it does not necessarily represent the speed of convergence to 
equilibrium.  In \cite{dss12}, where $q_x$ is set to be a common $q$ for all 
$x$, the average number of spin-flips per update is conceptually 
explained to be $O(|V|e^{-\beta q})$.  Here, we show an exact 
computation of the SCA transition probability and approximate it by a 
binomial expansion, from which we claim that the actual average number of 
spin-flips per update is much smaller than $O(|V|e^{-\beta q})$.

First, we recall equation  \refeq{PSCAdef}.  Notice that
\begin{align}
\frac{e^{\frac\beta2 (\tilde h_x(\sigmavec)+q_x\sigma_x)\tau_x}}{2\cosh(\frac
 \beta2(\tilde h_x(\sigmavec)+q_x\sigma_x))}
&=\frac{e^{-\frac\beta2(\tilde h_x(\sigmavec)\sigma_x+q_x)}\ind{x\in
 D_{\sigmavec,\tauvec}}}{2\cosh(\frac\beta2 (\tilde h_x(\sigmavec)+q_x
 \sigma_x))}+\frac{e^{\frac\beta2 (\tilde h_x(\sigmavec)\sigma_x+q_x)}\ind{x\in
 V\setminus D_{\sigmavec,\tauvec}}}{2\cosh(\frac\beta2(\tilde h_x(\sigmavec)
 +q_x\sigma_x))}.
\end{align}
Isolating the $\qvec$-dependence, we can rewrite the first term on the right-hand side as 
\begin{align}
\frac{e^{-\frac\beta2(\tilde h_x(\sigmavec)\sigma_x+q_x)}\ind{x\in D_{\sigmavec,
 \tauvec}}}{2\cosh(\frac\beta2 (\tilde h_x(\sigmavec)+q_x\sigma_x))}=
 \underbrace{\frac{e^{-\frac\beta2q_x}\cosh(\frac\beta2\tilde h_x(\sigmavec))}
 {\cosh(\frac\beta2(\tilde h_x(\sigmavec)+q_x\sigma_x))}}_{\equiv\,\vep_x
 (\sigmavec)}\,\underbrace{\frac{e^{-\frac\beta2\tilde h_x(\sigmavec) \sigma_x}}
 {2\cosh(\frac\beta2 \tilde h_x(\sigmavec))}}_{\equiv\,p_x(\sigmavec)}\ind{x\in
 D_{\sigmavec,\tauvec}},
\end{align}
and the second term as $(1-\vep_x(\sigmavec)p_x(\sigmavec))
\ind{x\in V\setminus D_{\sigmavec,\tauvec}}$.  As a result, we obtain
\begin{align}
\PSCA_{\beta,\qvec}(\sigmavec,\tauvec)=\prod_{x\in D_{\sigmavec,\tauvec}}\big(\vep_x
 (\sigmavec)p_x(\sigmavec)\big)\prod_{y\in V\setminus D_{\sigmavec,\tauvec}}
 \big(1-\vep_y(\sigmavec)p_y(\sigmavec)\big).
\end{align}
Suppose that $\vep_x(\sigmavec)$ is independent of $x$ and $\sigmavec$, which 
is of course untrue, and simply denote it by $\vep=O(e^{-\beta q})$.  Then we can 
rewrite $\PSCA_{\beta,\qvec}(\sigmavec,\tauvec)$ as
\begin{align}\lbeq{binom-exp}
\PSCA_{\beta,\qvec}(\sigmavec,\tauvec)
&\simeq\prod_{x\in D_{\sigmavec,\tauvec}}\big(
 \vep p_x(\sigmavec)\big)\prod_{y\in V\setminus D_{\sigmavec,\tauvec}}\Big((1
 -\vep)+\vep\big(1-p_y(\sigmavec)\big)\Big)\nn\\
&=\prod_{x\in D_{\sigmavec,\tauvec}}\big(\vep p_x(\sigmavec)\big)\sum_{S:
 D_{\sigmavec,\tauvec}\subset S\subset V} (1-\vep)^{|V\setminus S|}\prod_{y\in
 S\setminus D_{\sigmavec,\tauvec}}\Big(\vep\big(1-p_y(\sigmavec)\big)\Big)\nn\\
&=\sum_{S:D_{\sigmavec,\tauvec}\subset S\subset V}\vep^{|S|} (1-\vep)^{|V\setminus S|}
 \prod_{x\in D_{\sigmavec,\tauvec}}p_x(\sigmavec)\prod_{y\in S
 \setminus D_{\sigmavec,\tauvec}}\big(1-p_y(\sigmavec)\big).
\end{align}
This implies that the transition from $\sigmavec$ to $\tauvec$ can be seen as 
determining the binomial subset $D_{\sigmavec,\tauvec}\subset S\subset V$ with 
parameter $\vep$ and then changing each spin at $x\in D_{\sigmavec,\tauvec}$ 
with probability $p_x(\sigmavec)$.  Therefore, $|V|\vep$ is much larger 
than the actual average number of spin-flips per update.

Currently, the authors are investigating an MCMC inspired 
by \refeq{binom-exp} with a constant $\vep\in(0,1]$.  Some numerical results 
have shown better performance than Glauber dynamics and SCA in finding ground states for several 
problems. For more details, see Section \ref{sec:eSCA}. 
\item
In fact, we can compute the average number $E^*[|D_{\sigmavec,X}|]
\equiv\sum_{\tauvec}|D_{\sigmavec,\tauvec}|P^*(\sigmavec,\tauvec)$ of 
spin-flips per update, where $X$ is an $\Omega$-valued random variable whose 
law is $P^*(\sigmavec,\cdot)$.  For Glauber, we have
\begin{align}
\EG_\beta[|D_{\sigmavec,X}|]=\sum_{x\in V}\PG_\beta(\sigmavec,\sigmavec^x)
 =\frac1{|V|}\sum_{x\in V}\frac{e^{-\beta\tilde h_x(\sigmavec)
 \sigma_x}}{2\cosh(\beta\tilde h_x(\sigmavec))}
 =\frac1{|V|}\sum_{x\in V}\frac1{e^{2\beta\tilde h_x(\sigmavec)
 \sigma_x}+1}.
\end{align}
For the SCA, on the other hand, since 
$|D_{\sigmavec,\tauvec}|=\sum_{x\in V}\ind{\sigma_x\ne\tau_x}$, we have
\begin{align}
\ESCA_{\beta,\qvec}[|D_{\sigmavec,X}|]=\sum_{x\in V}\sum_{\tauvec:\tau_x\ne
 \sigma_x}\PSCA_{\beta,\qvec}(\sigmavec,\tauvec)&=\sum_{x\in V}\frac{e^{-\frac
 \beta2(\tilde h_x(\sigmavec)\sigma_x+q_x)}}{2\cosh(\frac\beta2(\tilde h_x
 (\sigmavec)\sigma_x+q_x))}\nn\\
&=\sum_{x\in V}\frac1{e^{\beta (\tilde h_x(\sigmavec)\sigma_x+q_x)}+1}.
\end{align}
Therefore, the SCA has more spin-flips per update than Glauber, if 
$|V|e^{2\beta\tilde h_x(\sigmavec)\sigma_x}\ge e^{\beta(\tilde h_x(\sigmavec)
\sigma_x+q_x)}$ for any $x\in V$ and $\sigmavec\in\Omega$, which is true 
when the temperature is sufficiently high such that
\begin{align}\lbeq{moreflipscond}
\max_{x\in V}\frac\beta2\bigg(q_x+|h_x|+\sum_{y\in V}|J_{x,y}|\bigg)\le\log
 \sqrt{|V|}.
\end{align}
Compare this with the condition \refeq{mixingcond}, which is independent of 
$\{h_x\}_{x\in V}$, hence better than \refeq{moreflipscond} in this respect.  
On the other hand, the bound in \refeq{moreflipscond} can be made large 
as $|V|$ increases, while it is always 1 in \refeq{mixingcond}.
\end{enumerate}
}
\end{rmk}

\section{Simulated annealing for the SCA}\label{sec:SCAannealing}
In this section, we show that, under a logarithmic cooling schedule 
$\beta_t\propto\log t$, the simulated annealing for the SCA weakly converges to 
the uniform distribution over GS.  To do so, we introduce the
Dobrushin's ergodic coefficient $\delta(P)$ of the transition matrix 
$[P(\sigmavec,\tauvec)]_{\Omega\times\Omega}$ as
\begin{align}\lbeq{Dobrushin}
\delta(P)=\max_{\sigmavec,\tauvec\in\Omega}\|P(\sigmavec,\cdot)-P(\tauvec,\cdot)
 \|_{\sss\rm TV}\equiv1-\min_{\sigmavec,\etavec}\sum_{\tauvec}P(\sigmavec,
 \tauvec)\wedge P(\etavec,\tauvec).
\end{align}

The following proposition is a summary of \cite[Theorems~6.8.2 \& 6.8.3]{b99}, 
but stated in our context.

\begin{shaded}
\begin{prp}\label{prp:ergodicity}
Let $\{X_n\}_{n=0}^\infty$ be a time-inhomogenous Markov chain on $\Omega$ 
generated by the transition probabilities $\{P_n\}_{n\in\N}$, i.e.,  
$P_n(\sigmavec,\tauvec)=\mP(X_n=\tauvec|X_{n-1}=\sigmavec)$.  Let 
$\{\pi_n\}_{n\in\N}$ be their respective equilibrium distributions, i.e., 
$\pi_n=\pi_nP_n$ for each $n\in\N$.  If 
\begin{align}\lbeq{weak->strong}
\sum_{n=1}^\infty\|\pi_{n+1}-\pi_n\|_{\sss\rm TV}<\infty,
\end{align}
and if there is a strictly increasing sequence $\{n_j\}_{j\in\N}\subset\N$ 
such that
\begin{align}\lbeq{weak-erg}
\sum_{j=1}^\infty\Big(1-\delta(P_{n_j}P_{n_j+1}\cdots P_{n_{j+1}-1})\Big)=\infty,
\end{align}
then there is a probability distribution $\pi$ on $\Omega$ such that, for any 
$j\in\N$,
\begin{align}\lbeq{strong-erg}
\lim_{n\uparrow\infty}\sup_\mu\|\mu P_j\cdots P_n-\pi\|_{\sss\rm TV}=0,
\end{align}
where the supremum is taken over the initial distribution on $\Omega$.
\end{prp}
\end{shaded}

The second assumption \refeq{weak-erg} is a necessary and sufficient condition 
for the Markov chain to be weakly ergodic \cite[Definition~6.8.1]{b99}: 
for any $j\in\N$,
\begin{align}
\lim_{n\uparrow\infty}\sup_{\mu,\nu}\|\mu P_j\cdots P_n-\nu P_j\cdots P_n
 \|_{\sss\rm TV}=0.
\end{align} 
On the other hand, if \refeq{strong-erg} holds, then the Markov chain is called 
strongly ergodic \cite[Definition~6.8.2]{b99}.  The first assumption 
\refeq{weak->strong} is to guarantee strong ergodicity from weak ergodicity, 
as well as the existence of the limiting measure $\pi$.

To apply this proposition to the SCA, it is crucial to find a cooling schedule 
$\{\beta_t\}_{t\in\N}$ under which the two assumptions 
\refeq{weak->strong}--\refeq{weak-erg} hold, and to show that the limiting 
measure is the uniform distribution $\piG_\infty$ over GS.  
Here is the main statement on the simulated annealing for the SCA.

\begin{shaded}
\begin{thm}\label{thm:SAforSCA}
Suppose that the pinning parameters $\qvec$ satisfy the condition 
\refeq{MA-pinning-parameter}.  For any non-decreasing sequence 
$\{\beta_t\}_{t\in\N}$ satisfying $\lim_{t\uparrow\infty}\beta_t=\infty$, 
we have
\begin{align}\lbeq{SCAweak->strong}
\sum_{t=1}^\infty\|\piSCA_{\beta_{t+1},\qvec}-\piSCA_{\beta_t,\qvec}\|_{\sss
 \rm TV}<\infty,&&
\lim_{t\uparrow\infty}\|\piSCA_{\beta_t,\qvec}-\piG_\infty\|_{\sss\rm TV}=0.
\end{align}
In particular, if we choose $\{\beta_t\}_{t\in\N}$ as
\begin{align}\lbeq{SCA-cooling-schedule}
\beta_t=\frac{\log t}\Gamma,&&
\Gamma=\sum_{x \in V} \Gamma_x,&&
\Gamma_x =q_x+|h_x|+\sum_{y\in V}|J_{x,y}|,
\end{align}
then we obtain
\begin{align}\lbeq{SCAweak-erg}
\sum_{t=1}^\infty\big(1-\delta(\PSCA_{\beta_t,\qvec})\big)=\infty.
\end{align}
As a result, for any initial $j\in\N$,
\begin{align}\label{SA:Sergodic}
\lim_{t\to\infty}\sup_\mu\big\|\mu\PSCA_{\beta_j,\qvec}\PSCA_{\beta_{j+1},
 \qvec}\cdots\PSCA_{\beta_t,\qvec}-\piG_\infty\big\|_{\sss\rm TV} = 0.
\end{align}
\end{thm} 
\end{shaded}

\begin{proof}
Since \eqref{SA:Sergodic} is an immediate consequence of 
Proposition~\ref{prp:ergodicity}, \refeq{SCAweak->strong} and 
\refeq{SCAweak-erg}, it remains to show \refeq{SCAweak->strong} and 
\refeq{SCAweak-erg}.

To show \refeq{SCAweak->strong}, we first define
\begin{align}\label{SA:mu}
\mu_\beta(\sigmavec,\tauvec)=\frac{e^{-\beta\tilde H(\sigmavec,\tauvec)}}
 {\sum_{\xivec,\etavec}e^{-\beta\tilde H(\xivec,\etavec)}}\equiv\frac{e^{-\beta
 (\tilde H(\sigmavec,\tauvec)-m)}}{\sum_{\xivec,\etavec}e^{-\beta(\tilde H
 (\xivec,\etavec) - m)}},
\end{align}
where $m=\min_{\sigmavec,\etavec}\tilde H(\sigmavec,\etavec)$.  
Since $\qvec$ is chosen to satisfy \refeq{mindiag}, we can conclude that
\begin{align}\label{SA:mu2}
\mu_\beta(\sigmavec,\tauvec)=\frac{e^{-\beta(\tilde H(\sigmavec,\tauvec)-m)}}
 {|\mathrm{GS}|+\sum_{\xivec,\etavec:\tilde H(\xivec,\etavec)>m}e^{-\beta(\tilde
 H(\xivec,\etavec) - m)}}\xrightarrow[\beta\uparrow\infty]{}\underbrace{\frac{
 \ind{\sigmavec\in\mathrm{GS}}}{|\mathrm{GS}|}}_{\piG_\infty(\sigmavec)}
 \delta_{\sigmavec,\tauvec}.
\end{align}
Summing this over $\tauvec\in\Omega\equiv\{\pm1\}^V$ yields the second 
relation in \refeq{SCAweak->strong}.  To show the first relation in 
\refeq{SCAweak->strong},  we note that
\begin{align}\label{partderiv}
\frac{\partial\mu_\beta(\sigmavec,\tauvec)}{\partial\beta}=\Big(\mE_{\mu_\beta}
 [\tilde H]-\tilde H(\sigmavec, \tauvec)\Big)\mu_\beta(\sigmavec,\tauvec),
\end{align}
and that $\mE_{\mu_\beta}[\tilde H]\equiv\sum_{\sigmavec,\tauvec}\tilde H
(\sigmavec,\tauvec)\,\mu_\beta(\sigmavec,\tauvec)$ tends to $m$ as 
$\beta\uparrow\infty$, due to \eqref{SA:mu2}.  Therefore, 
$\frac\partial{\partial\beta}\mu_\beta(\sigmavec,\tauvec)>0$ for all $\beta$ 
if $\tilde H(\sigmavec,\tauvec)=m$, while it is negative for sufficiently large 
$\beta$ if $\tilde H(\sigmavec,\tauvec)>m$.  Let $n\in\N$ be such that, 
as long as $\beta\ge\beta_n$, \eqref{partderiv} is negative for all pairs 
$(\sigmavec,\tauvec)$ satisfying $\tilde H(\sigmavec,\tauvec)>m$.  As a result, 
\begin{align}
&\sum_{t=n}^N\|\piSCA_{\beta_{t+1},\qvec}-\piSCA_{\beta_t,\qvec}\|_\mathrm{TV}\nn\\
&=\frac12\sum_{\sigmavec \in \mathrm{GS}}\sum_{t =n}^N|\piSCA_{\beta_{t+1},\qvec}
 (\sigmavec)-\piSCA_{\beta_t,\qvec}(\sigmavec)|+\frac12\sum_{\sigmavec \notin
 \mathrm{GS}}\sum_{t =n}^N|\piSCA_{\beta_{t+1},\qvec}(\sigmavec)-\piSCA_{\beta_t,
 \qvec}(\sigmavec)|\nn\\
&\le\frac12\sum_{\sigmavec\in\mathrm{GS}}\sum_{t =n}^N\big(\mu_{\beta_{t+1}}
 (\sigmavec,\sigmavec)-\mu_{\beta_t}(\sigmavec,\sigmavec)\big)+\frac12
 \sum_{\sigmavec\in\mathrm{GS}}\sum_{\tauvec\ne\sigmavec}\sum_{t=n}^N\big(
 \mu_{\beta_t}(\sigmavec,\tauvec)-\mu_{\beta_{t+1}}(\sigmavec,\tauvec)\big)\nn\\
&\quad+\frac12\sum_{\sigmavec\notin\mathrm{GS}}\sum_{t=n}^N\big(\piSCA_{\beta_t,
 \qvec}(\sigmavec)-\piSCA_{\beta_{t+1},\qvec}(\sigmavec)\big)\nn\\
&=\frac12\sum_{\sigmavec\in\mathrm{GS}}\big(\mu_{\beta_{N+1}}(\sigmavec,\sigmavec)
 -\mu_{\beta_n}(\sigmavec,\sigmavec)\big)+\frac12\sum_{\sigmavec\in\mathrm{GS}}
 \sum_{\tauvec\ne\sigmavec}\big(\mu_{\beta_n}(\sigmavec,\tauvec)-\mu_{\beta_{N
 +1}}(\sigmavec,\tauvec)\big)\nn\\
&\quad+\frac12\sum_{\sigmavec\notin\mathrm{GS}}\big(\piSCA_{\beta_n,\qvec}
 (\sigmavec)-\piSCA_{\beta_{N+1},\qvec}(\sigmavec)\big)\nn\\
&\le\frac32
\end{align}
holds uniformly for $N\ge n$.  This completes the proof of 
\refeq{SCAweak->strong}.

To show \refeq{SCAweak-erg}, we use the following bound on 
$\PSCA_{\beta,\qvec}$, which holds uniformly in $(\sigmavec,\tauvec)$:
\begin{eqnarray}
\PSCA_{\beta,\qvec}(\sigmavec,\tauvec)\stackrel{\text{\refeq{PSCAdef}}}=\prod_{x
 \in V}\frac{e^{\frac\beta2(\tilde h_x(\sigmavec)+q_x\sigma_x)\tau_x}}{2\cosh
 (\frac\beta2(\tilde h_x(\sigmavec)+q_x\sigma_x))}&\ge&\prod_{x\in V}\frac1{1
 +e^{\beta|\tilde h_x(\sigmavec)+q_x\sigma_x|}}\nn\\
&\stackrel{\text{\refeq{SCA-cooling-schedule}}}\ge&\prod_{x\in V}\frac{e^{-\beta
 \Gamma_x}}2=\frac{e^{-\beta\Gamma}}{2^{|V|}}.
\end{eqnarray}
Then, by \refeq{Dobrushin}, we obtain
\begin{align}
\sum_{t=1}^\infty\big(1-\delta(\PSCA_{\beta_t,\qvec})\big)=\sum_{t=1}^\infty
 \min_{\sigmavec,\etavec}\sum_{\tauvec}\PSCA_{\beta_t,\qvec}(\sigmavec,\tauvec)
 \wedge\PSCA_{\beta_t,\qvec}(\etavec,\tauvec)\ge\sum_{t=1}^\infty e^{-\beta_t
 \Gamma},
\end{align}
which diverges, as required, under the cooling schedule 
\refeq{SCA-cooling-schedule}.  This completes the proof of the theorem.
\end{proof}

\begin{rmk}\label{rmk:SAforSCA}
{\rm
\begin{enumerate}[(i)]
\item
{The main message contained in the above theorem is that, in order to achieve weak 
convergence to the uniform distribution over the ground states, it is enough for 
the temperature to drop no faster than $1/\log t$ with a large multiplicative 
constant $\Gamma$.  For logarithmic schedules, due to our approach, it is not trivial to ensure whether the value of  $\Gamma$ as in
\refeq{SCA-cooling-schedule} is optimal for the SCA, contrasting with the Metropolis dynamics, whose optimal value can be  theoretically determined, see 
\cite{h88}.}
\item
Simulated annealing with a logarithmic cooling schedule may not be so 
practical in finding a ground state within a feasible amount of time.  Instead, 
an exponential cooling schedule is often used in engineering.  In 
\cite{IEEE2023,STATICA}, we have developed annealing processors called Amorphica and STATICA, 
that rely on the SCA with an exponential schedule.  Experimental results have shown 
faster in searching for a ground state than conventional simulated annealing (based on the Glauber 
dynamics with an exponential schedule)
%, just as mentioned in Remark~\ref{rmk:mixing}(i), 
and better performance in finding solutions to 
a max-cut problem.

The authors are investigating the use of 
exponential schedules.  We do not expect weak convergence to the uniform 
distribution over the ground states.  Instead, we want to evaluate the 
probability that the SCA with an exponential schedule reaches a spin 
configuration $\sigmavec$ such that $H(\sigmavec)-\min H$ is within a given 
error margin.  
A similar problem was considered by Catoni \cite{c92} for the Metropolis 
dynamics with an exponential schedule.
\item
However, this may imply that we have not yet been able to make the most of 
the SCA's independent multi-spin flip rule for better cooling schedules.  
The use of the total-variation distance may be one of the reasons why we have 
to impose such tight conditions on the temperature; if two measures $\mu$ and 
$\nu$ on $\Omega$ are close in total variation, then 
$|\mu(\sigmavec)-\nu(\sigmavec)|$ must be small at every $\sigmavec\in\Omega$.  
We should keep in mind that the most important thing in combinatorial optmization 
is to know the ordering among spin configurations, and not to perfectly fit $\piG_\beta$ by 
$\piSCA_{\beta,\qvec}$.  For example, $\piSCA_{\beta,\qvec}$ does not have to 
be close to $\piG_\beta$ in total variation, as long as we can say instead that 
$H(\sigmavec)\le H(\tauvec)$ (or equivalently 
$\piG_\beta(\sigmavec)\ge\piG_\beta(\tauvec)$) whenever 
$\piSCA_{\beta,\qvec}(\sigmavec)\ge\piSCA_{\beta,\qvec}(\tauvec)$ 
(see Figure~\ref{fig:dist}).
\begin{figure}
\begin{align*}
\includegraphics[scale=.35]{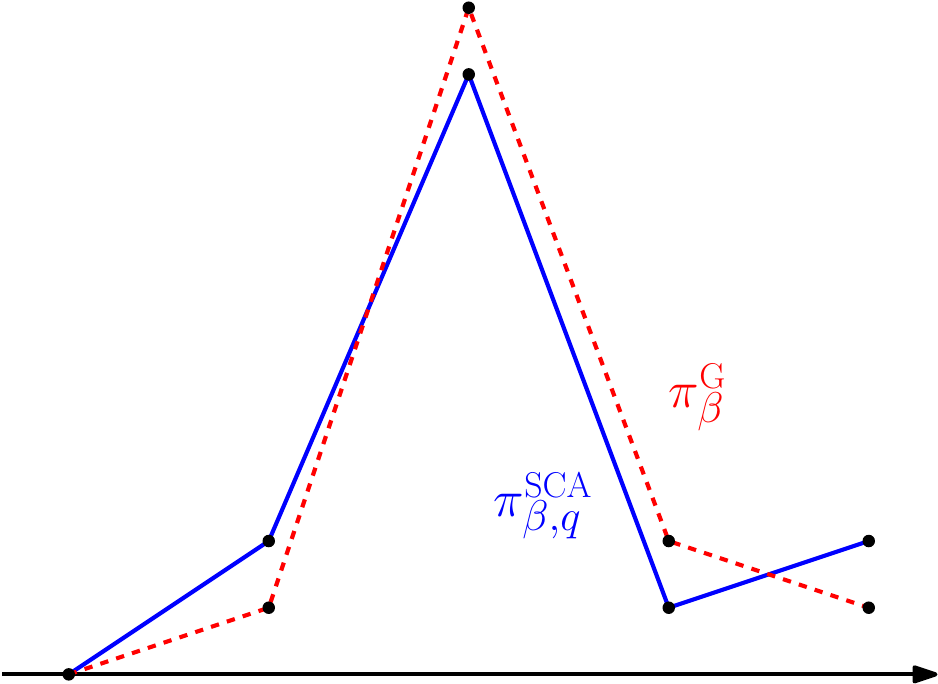}&&&&
\includegraphics[scale=.35]{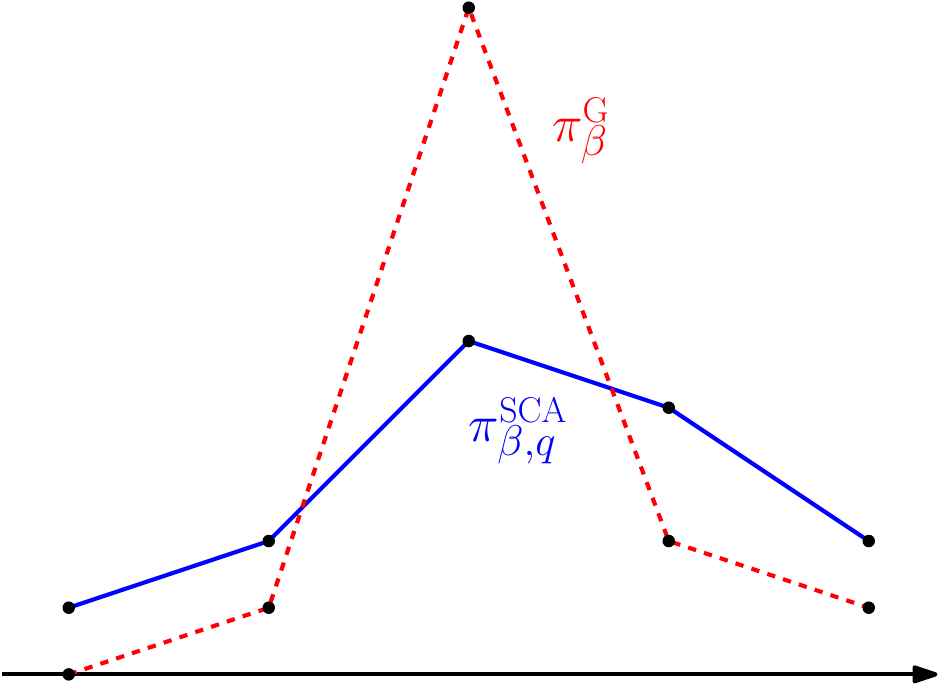}
\end{align*}
\caption{
On the left, $\|\piSCA_{\beta,\qvec}-\piG_\beta\|_{\sss\rm TV}$ is small, 
but the ordering among spin configurations is not preserved.  
On the right, $\|\piSCA_{\beta,\qvec}-\piG_\beta\|_{\sss\rm TV}$ is not 
small, but the ordering among spin configurations is preserved.}\label{fig:dist}
\end{figure}

In \cite{hkks19}, we introduced a slightly relaxed version of this notion of 
closeness, which is also used in the stability analysis \cite{fstu21}.  
Given an error ratio $\vep\in(0,1)$, the SCA equilibrium measure 
$\piSCA_{\beta,\qvec}$ is said to be $\vep$-close to the target Gibbs 
$\piG_\beta$ in the sense of order-preservation if
\begin{align}\lbeq{closedef}
\piSCA_{\beta,\qvec}(\sigmavec)\ge\piSCA_{\beta,\qvec}(\tauvec)
 \quad\Rightarrow\quad H(\sigmavec)\le H(\tauvec)+\vep R_H~~
 \Big(\Leftrightarrow~\piG_\beta(\sigmavec)\ge\piG_\beta(\tauvec)e^{-\beta
 \vep R_H}\Big),
\end{align}
where $R_H\equiv\max_{\sigmavec,\tauvec}|H(\sigmavec)-H(\tauvec)|$ is the 
range of the Hamiltonian.
By simple arithmetic \cite{hkks19} (with a little care needed due to the 
difference in the definition of $\tilde H$), we can show that 
$\piSCA_{\beta,\qvec}$ is $\vep$-close to $\piG_\beta$ if, for all $x\in V$,
\begin{align}\lbeq{lowerbd}
q_x\ge|h_x|+\sum_y|J_{x,y}|+\frac1\beta\log\frac{2|V|(|h_x|+\sum_y|J_{x,y}|)}
 {\vep R_H}.
\end{align}
Unfortunately, this is not better than \refeq{MA-pinning-parameter}, which 
we recall is a sufficient condition for $\piSCA_{\beta,\qvec}$ to attain the 
highest peaks over GS, and not anywhere else.  Since 
$|V|(|h_x|+\sum_y|J_{x,y}|)/R_H$ in the logarithmic term in 
\refeq{lowerbd} is presumably of order 1, we can say that, if the assumption 
\refeq{MA-pinning-parameter} is slightly tightened to 
$q_x\ge|h_x|+\sum_y|J_{x,y}|+O_\vep(\beta^{-1})$, then the SCA 
can be used to find not only the best options in combinatorial optimization, 
but also the second- and third-best options, etc.  In an ongoing project, we are 
also aiming to improve the cooling schedule under the new notion of 
closeness.
\end{enumerate}
}
\end{rmk}

\section{Comparisons and simulations}\label{sec:eSCA}

Based on the discussion from Remark \ref{rmk:mixing}(ii), let us propose a new algorithm derived from the SCA 
studied in this paper and make a quick comparison  regarding their effectiveness in obtaining the ground states. Given a fixed inverse temperature $\beta \geq 0$ and a number $\vep \in [0,1]$, let the 
transition matrix of the $\vep$-SCA be defined by
\begin{equation}
P_{\beta,\vep}(\sigmavec,\tauvec) = \prod_{x\in D_{\sigmavec,\tauvec}}\big(\vep p_x(\sigmavec)\big)\prod_{y\in V\setminus D_{\sigmavec,\tauvec}}
 \big(1-\vep p_y(\sigmavec)\big),
\end{equation}
where we recall that
\begin{equation}\lbeq{flipprobability}
p_x(\sigmavec) = \frac{e^{-\frac\beta2\tilde h_x(\sigmavec) \sigma_x}}
 {2\cosh(\frac\beta2 \tilde h_x(\sigmavec))}
\end{equation} 
is the probability of flipping the spin $\sigma_x$ from the configuration $\sigmavec$ disregarding a pinning parameter at $x$. Note that $1-\vep p_x(\sigmavec) = (1 - \vep) + \vep (1 - p_x(\sigmavec))$. Therefore, we can visualize this new algorithm by decomposing it into two steps: in the first step, the spins which are eligible to be flipped are selected independently at random, where each spin is selected with probability $\vep$, while it remains unchanged with probability $1 - \vep$; in the second step all spins which were selected in the previous step are updated simultaneously and independently, where the probability of flipping the spin at $x$ is equal to $p_x(\sigmavec)$. 

Note that,  in the particular case where $\vep = 1$, the update rule we have just introduced coincides with the SCA transition probability without the pinning parameters. Our experience \cite{SSS21,TISCIE23,IEEE2023} has shown that, for a fixed Hamiltonian, exponential cooling schedule and number of Markov chain steps, the simulated annealing algorithm based on $\vep$-SCA  with appropriately chosen parameter $\vep$  surpasses the performance of the algorithms based on SCA and Glauber dynamics with respect to the success probability in obtaining an approximation for a ground state. Later in this section, we will return to the question of how the performance of this algorithm is affected by the value of the parameter $\vep$.  

\begin{figure}[t!]
        \centering
        \begin{subfigure}[b]{0.49\textwidth}
                \includegraphics[width=\textwidth]{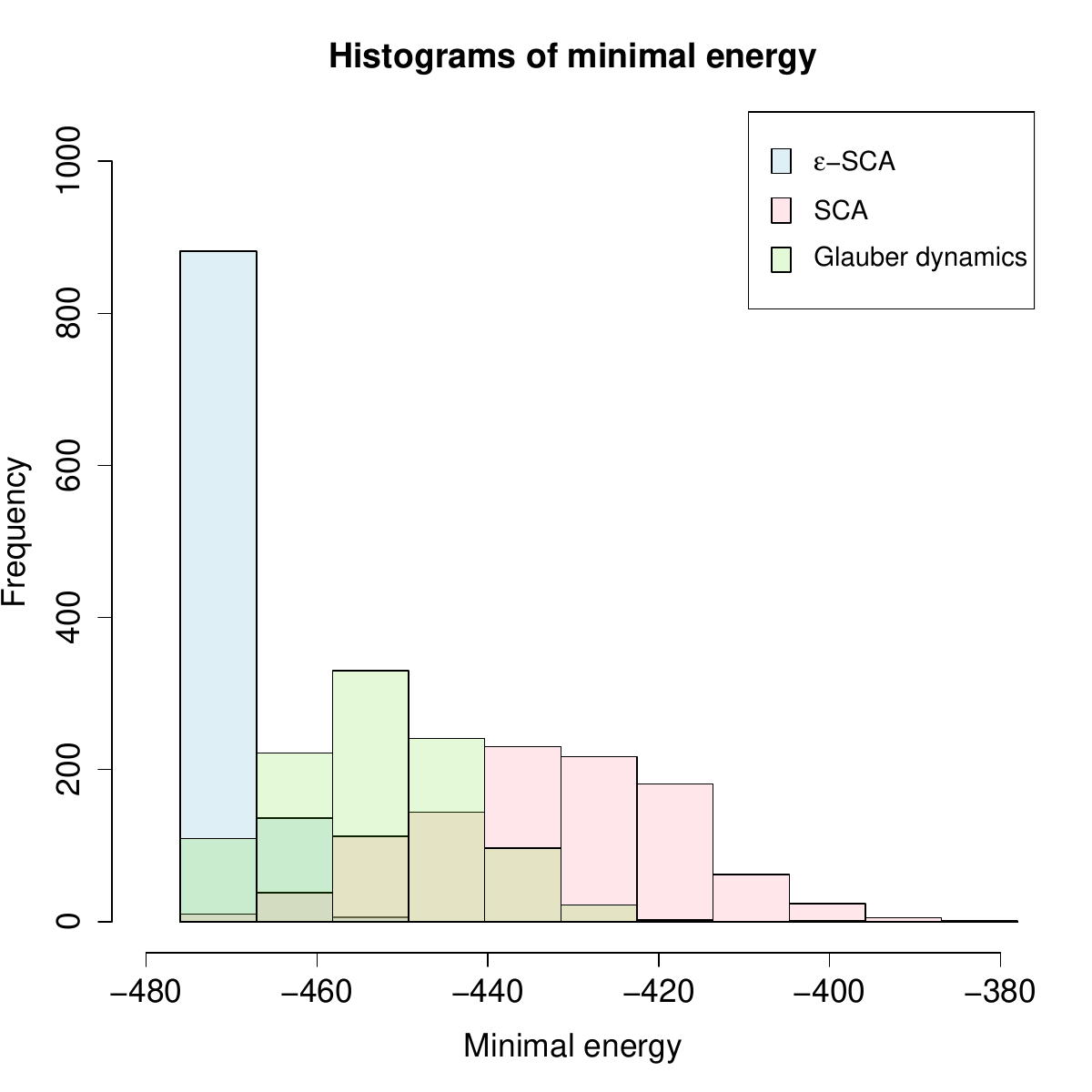}
                \caption{Max-cut problem}
                \label{fig:MCeSCA2}
        \end{subfigure} 
        \begin{subfigure}[b]{0.49\textwidth}
                \includegraphics[width=\textwidth]{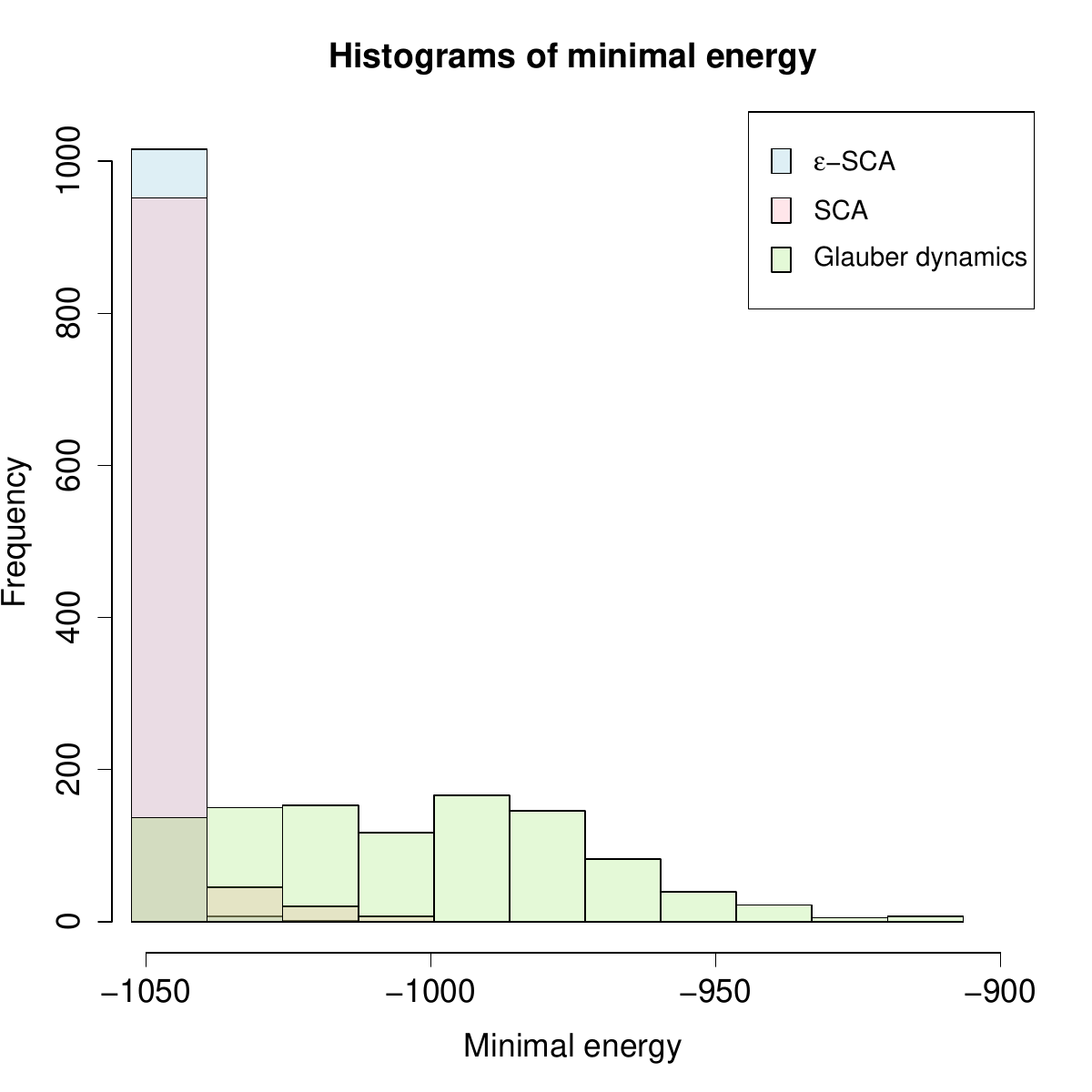}
                \caption{Spin-glass}
                \label{fig:MCeSCA}
        \end{subfigure}       
      
    \caption{The histograms comparing the success rates of obtaining a low-energy configuration for the three algorithms.}\label{fig:SIM}
\end{figure}

Now, let us make a comparison between the performances of the simulated annealing algorithms based on Glauber dynamics, SCA (with pinning parameters uniformly taken as $q_x = \frac{\lambda}{2}$) and $\vep$-SCA applied to the particular problems of determining the maximum cut of a given graph and to the minimization of a spin-glass Hamiltonian. Even though we only have rigorous results that justify the application of logarithmic cooling schedules to the Glauber dynamics \cite{AK1989,b99} and to the SCA, our practice has shown that exponential cooling schedules may also work for all these three algorithms, but we still do not have a solid theoretical justification for that. Let us restrict ourselves to the case where we have $N = 128$ spin variables. The plots from Figure \ref{fig:SIM} illustrate the histograms of the minimal energies achieved by running $M = 1024$ independent trials starting from randomly chosen initial spin configurations, where for each trial, we applied  $L = 20000$ Markov chain steps and considered the exponential cooling schedule with initial temperature  $T_\text{init} = 1000$ and final temperature $T_\text{fin} = 0.05$, explicitly, we considered
\begin{equation}\label{eq:cooling}
\frac{1}{\beta_t} = T_\text{init} \left(\frac{T_\text{fin}}{T_\text{init}}\right)^{\frac{t-1}{L-1}}
\end{equation}
for $t = 1, 2, \dots, L$.  
%%%%%%%%%%%%%%%
\begin{table}[b!]
\caption{The hitting rate $\rho_\text{min}$ to a configuration with the smallest energy among all those obtained by the three algorithms.}
\label{tab:perf1}
\centering
\begin{tabular}{cccccccc}
\hline
\multirow{2}{*}{Model}                          & \multicolumn{7}{c}{$\rho_\text{min}$} \\ \cline{2-8} 
                                                &   & Glauber  &  & SCA   &  & $\vep$-SCA &  \\ \hline \hline
Max-cut &   & 4.69\%      &  & 0\%    &  & 83.50\%     &  \\
Spin-glass                                      &   & 3.32\%     &  & 38.67\% &  & 61.33\%      &  \\ \hline
\end{tabular}
\end{table}
%%%%%%%%%%%%%%%%%%%
First, we fixed a randomly generated  Erd\" os-R\' enyi random graph $G(N,p)$ with $N = 128$ vertices and edge probability $p = 0.25$, and considered the Hamiltonian corresponding to the max-cut problem on $G(N,p)$, where $h_x = 0$ for every vertex $x$, and  $J_{x,y} = -1$ if $\{x,y\}$ is an edge of the graph, and $J_{x,y} =0$ otherwise.  The smallest energy obtained by the Glauber dynamics and the $\vep$-SCA with parameter $\vep = 0.3$ was equal $-476$, reached with success rates of $4.69\%$ and $83.50\%$, respectively, while the smallest energy obtained by the SCA was $-474$, reached with a success rate of $0.19\%$. Later, we fixed a spin-glass Hamiltonian on the complete graph $K_N$, where $N = 128$, whose spin-spin coupling constants were taken as the realizations of mutually independent standard normal random variables. In this case, the three algorithms obtained the same lowest energy, equal to
$-1052.57$, with a success rate of $3.32\%$ for the Glauber dynamics, $38.67\%$ for the SCA, while the obtained for the $\vep$-SCA with parameter $\vep = 0.8$ was of $61.33\%$. All such results are summarized in Table \ref{tab:perf1}. These  examples are consistent with our observations in \cite{SSS21,TISCIE23,IEEE2023}, where it was possible to conclude that the SCA outperforms the Glauber dynamics in certain scenarios where anti-ferromagnetic interactions are not prevalent, and the $\vep$-SCA outperforms both algorithms in all the studied scenarios.

\begin{figure}[t]
\centering
        \begin{subfigure}[b]{0.49\textwidth}
                \includegraphics[width=\textwidth]{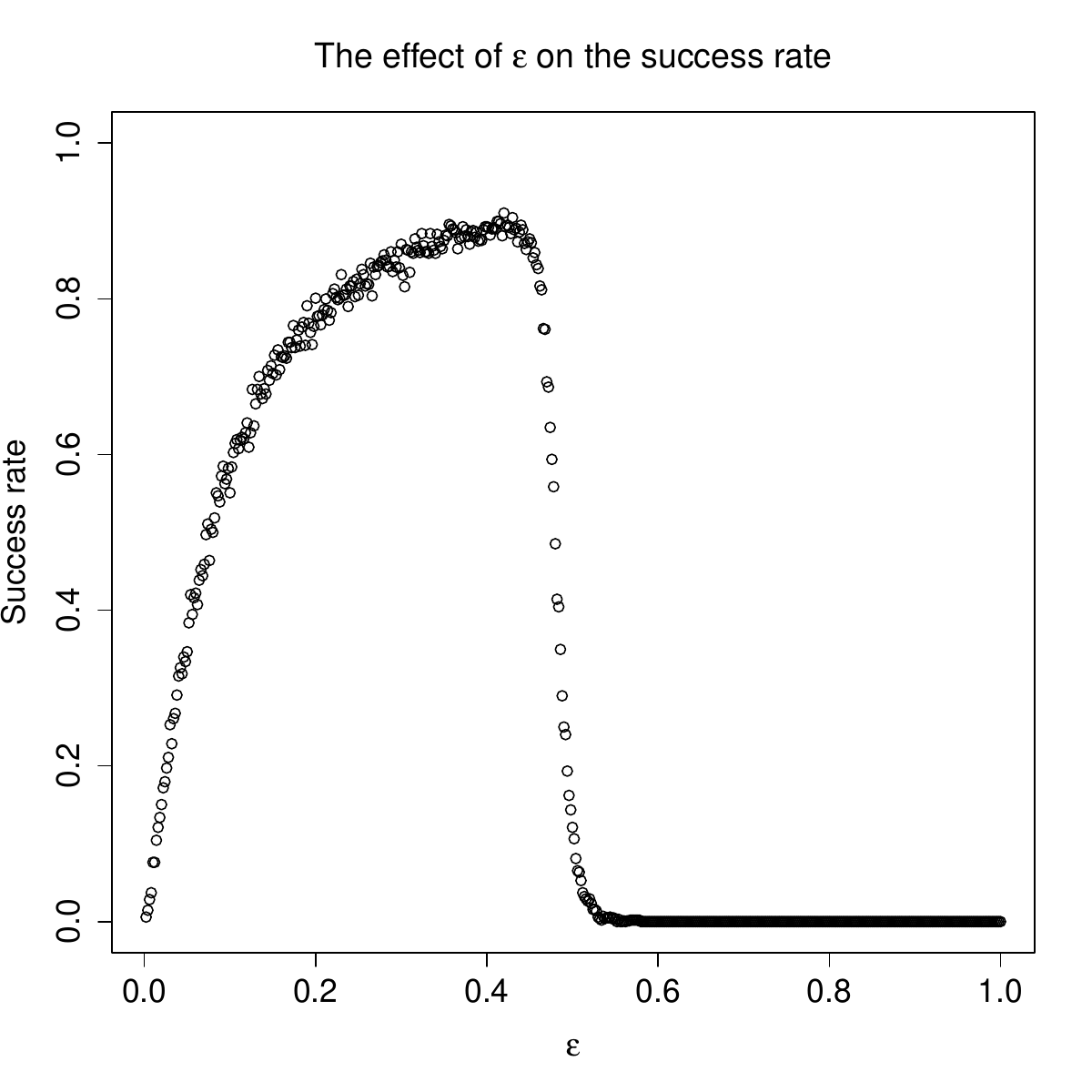}
                \caption{Max-cut problem}
                \label{fig:SRMC}
        \end{subfigure}
        \begin{subfigure}[b]{0.49\textwidth}
                \includegraphics[width=\textwidth]{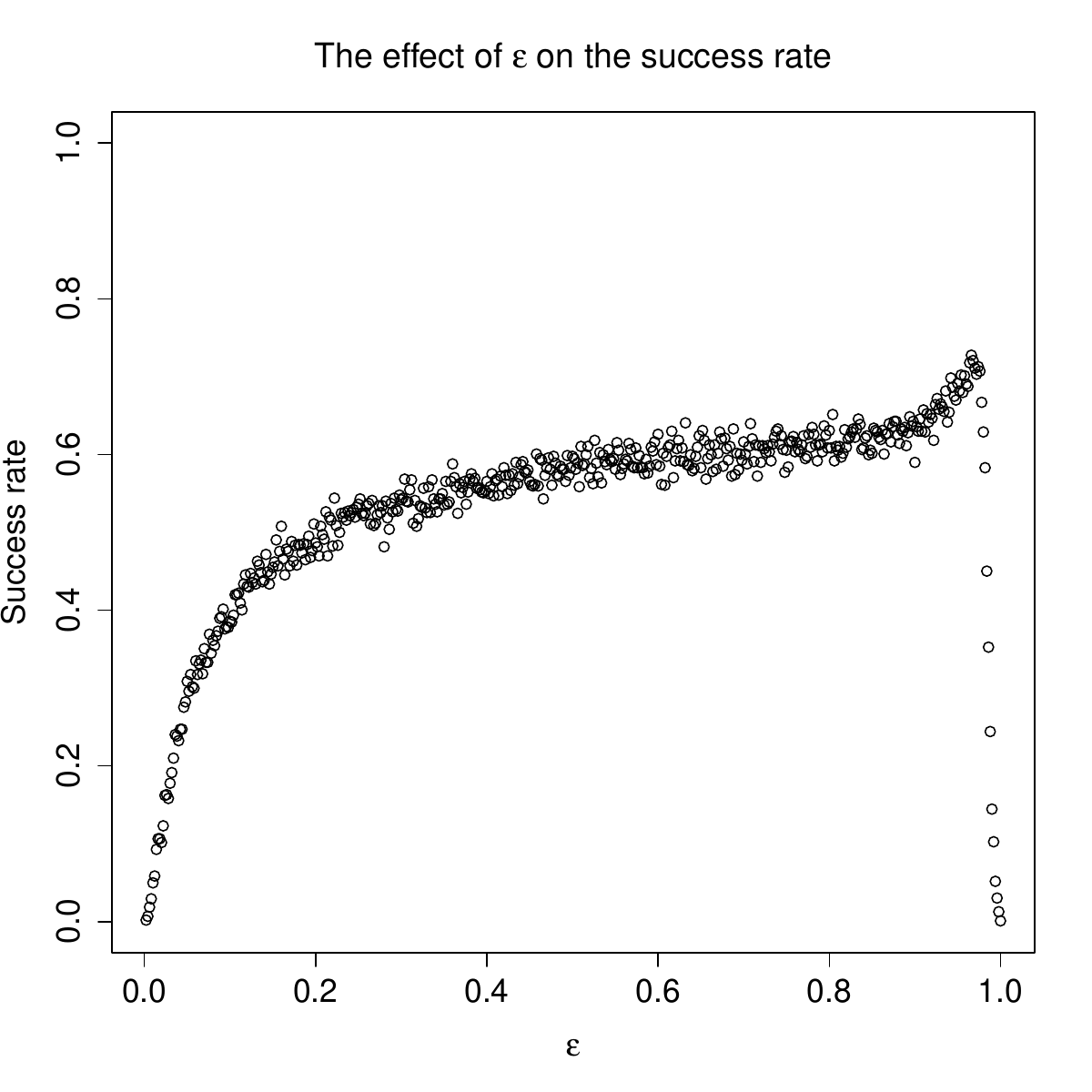}
                \caption{Spin-glass Hamiltonian}
                \label{fig:SRSG}
        \end{subfigure}       
\caption{Success rate of the $\vep$-SCA algorithm in obtaining a ground state as a function of the parameter $\vep$.}
\end{figure}

 The role played by the parameter $\vep$ in the $\vep$-SCA is analogous to the one played by the pinning parameters $\qvec=\{q_x\}_{x\in V}$ in the SCA. However, the difference is that the pinning effect for the SCA gets stronger as we decrease the temperature, so, the system will tend to flip less and less spins and might get stuck in an energetic local minimum. Regarding the $\vep$-SCA, due to the absence of the pinning parameters in the local transition probabilities and due to the effect of $\vep$ not being influenced by the temperature, the probability of flipping a certain spin will be larger compared to the SCA. Thus, this new algorithm allows the system to visit more configurations, especially at low temperatures, while preventing it from getting stuck in a local minimum. Nevertheless, a rigorous explanation that indicates what leads the system to converge to a ground state is still under investigation.

In order to get some intuition about the dependence on $\vep$ of the success rate of reaching a ground state, we considered the same max-cut problem and spin-glass Hamiltonian again and performed $M=1024$ $\vep$-SCA annealing trials for different values of $\vep$ with the same cooling schedule and number of Markov chain steps  as before. Typically, for Hamiltonians containing predominantly anti-ferromagnetic pairwise interactions (such as the Hamiltonian corresponding to the max-cut problem), such parameter $\vep$ has to be taken relatively small since the system tends to show an oscillatory behavior and does not converge to a ground state as we allow a larger number of spins to be flipped at a time, see Figure \ref{fig:SRMC}. On the other hand, for the spin-glass Hamiltonian, there is a tendency of growth of the success rate as the parameter $\vep$ increases. However, when $\vep$ gets sufficiently close to $1$, the algorithm behaves similarly to the SCA with pinning parameters $q_x = 0$, so the success rate decreases since the system will not necessarily converge to a ground state of the Hamiltonian, see Figure \ref{fig:SRSG}.
 Differently from the SCA, in which  we have a sufficient condition on the values for the pinning parameters that guarantee the convergence of the algorithm, it is still necessary to derive an analogous condition on $\vep$.  

 \begin{table}[b!]
\centering
\caption{The total annealing time $t_A$, measured in milliseconds, divided by the number of trials M corresponding to different hardware, algorithms, and number of vertices $N$.}
\label{tab:annealing}
\resizebox{\textwidth}{!}{%
\begin{tabular}{ccccccccccc} 
\hline
\multirow{3}{*}{M} & \multirow{3}{*}{Hardware} & \multirow{3}{*}{Algorithm} &  & \multicolumn{6}{c}{Total annealing time per replica [ms]}   &   \\ 
\cline{4-11}
                   &                           &                            &  & \multicolumn{6}{c}{N}                                       &   \\ 
\cline{5-10}
                   &                           &                            &  & 128    & 256    & 512     & 1024     & 2048     & 4096      &   \\ 
\hline\hline
1                  & CPU                       & Glauber                    &  & 60.80  & 235.66 & 898.11  & 3660.85  & 15815.50 & 71992.30  &   \\
1                  & CPU                       & SCA                        &  & 196.77 & 693.83 & 2533.02 & 9471.94  & 37684.70 & 144128.00 &   \\
1                  & CPU                       & $\vep$-SCA                 &  & 247.54 & 909.35 & 3793.08 & 13570.30 & 58396.30 & 232225.00 &   \\
1                  & GPU                       & Glauber                    &  & 131.15 & 139.82 & 152.62  & 189.56   & 777.99   & 2550.75   &   \\
1                  & GPU                       & SCA                        &  & 132.00 & 142.38 & 156.35  & 194.59   & 793.21   & 2562.43   &   \\
1                  & GPU                       & $\vep$-SCA                 &  & 132.08 & 141.42 & 155.31  & 194.33   & 794.06   & 2563.39   &   \\
128                & CPU                       & Glauber                    &  & 59.88  & 234.59 & 914.16  & 3641.63  & 15928.00 & 71687.70  &   \\
128                & CPU                       & SCA                        &  & 194.49 & 692.82 & 2531.55 & 9462.83  & 37452.90 & 143932.00 &   \\
128                & CPU                       & $\vep$-SCA                 &  & 241.92 & 914.47 & 3646.10 & 14401.20 & 59392.70 & 242698.00 &   \\
128                & GPU                       & Glauber                    &  & 2.40   & 4.77   & 14.10   & 47.77    & 577.06   & 2296.14   &   \\
128                & GPU                       & SCA                        &  & 2.70   & 5.48   & 14.98   & 49.04    & 578.86   & 2299.03   &   \\
128                & GPU                       & $\vep$-SCA                 &  & 2.67   & 5.49   & 14.97   & 48.44    & 578.86   & 2299.03   &   \\
\hline
\end{tabular}
}
\end{table}

% \usepackage{multirow}
\begin{sidewaystable}
\centering
\caption{Comparison of the performances of simulated annealing based on Glauber dynamics, SCA, and $\vep$-SCA (with $\vep = 0.8$).}
\label{tab:comparison}
\resizebox{\textwidth}{!}{%
\begin{tabular}{cccccccccccccc}
\hline
\multirow{2}{*}{$N$} &
  \multirow{2}{*}{$L$} &
   &
  \multicolumn{3}{c}{Glauber dynamics} &
   &
  \multicolumn{3}{c}{SCA} &
   &
  \multicolumn{3}{c}{$\vep$-SCA} \\ \cline{4-6} \cline{8-10} \cline{12-14} 
 &
   &
   &
  $\rho_\text{min}$ (\%) &
  $H_{\text{mean}}^{\text{(min)}}$ &
  $t_A / M$ {[}ms{]} &
   &
  $\rho_\text{min}$ (\%) &
  $H_{\text{mean}}^{\text{(min)}}$ &
  $t_A / M$ {[}ms{]} &
   &
  $\rho_\text{min}$ (\%) &
  $H_{\text{mean}}^{\text{(min)}}$ &
  $t_A / M$ {[}ms{]} \\ \hline \hline
$128$  & $10^4$            &  & 1.07  & -985.32   & 0.73   &  & 21.19 & -1045.35  & 0.80   &  & 55.76 & -1051.50  & 0.81   \\
$128$  & $2.5 \times 10^4$ &  & 3.71  & -1011.35  & 1.83   &  & 43.75 & -1049.59  & 2.02   &  & 64.36 & -1051.85  & 2.02   \\
$128$  & $5 \times 10^4$   &  & 11.04 & -1024.97  & 3.66   &  & 54.30 & -1051.31  & 4.04   &  & 73.44 & -1052.11  & 4.05   \\
$128$  & $10^5$            &  & 22.75 & -1036.76  & 7.34   &  & 64.84 & -1051.88  & 8.10   &  & 84.86 & -1052.31  & 8.11   \\ \hline
$256$  & $10^4$            &  & 0     & -2767.51  & 1.94   &  & 9.67  & -3020.12  & 2.10   &  & 88.18 & -3040.03  & 2.10   \\
$256$  & $2.5 \times 10^4$ &  & 0.20  & -2843.18  & 4.85   &  & 21.39 & -3033.72  & 5.26   &  & 96.29 & -3041.69  & 5.24   \\
$256$  & $5 \times 10^4$   &  & 0.49  & -2892.30  & 9.69   &  & 42.09 & -3038.17  & 10.50  &  & 98.34 & -3041.91  & 10.49  \\
$256$  & $10^5$            &  & 2.05  & -2940.48  & 19.44  &  & 64.06 & -3040.45  & 21.06  &  & 98.93 & -3042.06  & 21.01  \\ \hline
$512$  & $10^4$            &  & 0     & -7793.74  & 6.62   &  & 0     & -8660.45  & 7.02   &  & 14.75 & -8746.99  & 6.98   \\
$512$  & $2.5 \times 10^4$ &  & 0     & -8073.91  & 16.55  &  & 0.10  & -8713.27  & 17.53  &  & 20.61 & -8761.51  & 17.50  \\
$512$  & $5 \times 10^4$   &  & 0     & -8204.84  & 33.12  &  & 1.17  & -8738.57  & 35.13  &  & 24.90 & -8765.14  & 35.07  \\
$512$  & $10^5$            &  & 0     & -8368.72  & 66.21  &  & 1.66  & -8753.09  & 70.22  &  & 34.18 & -8767.84  & 70.01  \\ \hline
$1024$ & $10^4$            &  & 0     & -21028.50 & 24.52  &  & 0     & -24362.40 & 24.98  &  & 0.49  & -24636.00 & 24.89  \\
$1024$ & $2.5 \times 10^4$ &  & 0     & -22276.60 & 61.38  &  & 0     & -24514.30 & 62.44  &  & 0.68  & -24688.20 & 62.27  \\
$1024$ & $5 \times 10^4$   &  & 0     & -22815.30 & 122.86 &  & 0     & -24585.60 & 125.01 &  & 0.59  & -24710.80 & 124.46 \\
$1024$ & $10^5$            &  & 0     & -23235.10 & 242.43 &  & 0     & -24635.70 & 247.16 &  & 0.29  & -24724.70 & 246.22 \\ \hline
\end{tabular}
}
\end{sidewaystable}

{As we mentioned previously in this paper, the parallel spin update, which is the most notorious feature shared by the SCA and the $\vep$-SCA, leads us to consider hardware accelerators that can fully take advantage of such parallel nature and be very decisive in solving large-scale combinatorial optimization problems within a shorter amount of time as compared to these algorithms implemented on an average consumer-level device. 
STATICA \cite{STATICA} and Amorphica \cite{IEEE2023} are cutting-edge annealing processors developed to accommodate SCA and $\vep$-SCA annealing algorithms and have achieved results with higher precision, energy efficiency, and shorter execution times, compared to other state-of-art devices. In order to test the efficiency of our algorithms implemented in consumer-level processors, we evaluated the annealing time of Glauber dynamics, SCA, and $\vep$-SCA, by executing each algorithm in C++ (for running on a CPU) and C++/CUDA (for running on a GPU). In this evaluation, we ran the C++ code on an Intel Core i7-9700 Processor and the GPU kernel on an NVIDIA GeForce RTX 2080 Ti.
Let $N$ be the number of spin variables involved and $M$ the number of annealing trials/replicas. The execution time on a CPU is approximately proportional to $N^2 \times M$. In contrast, thanks to the parallel computing capability, the execution efficiency can be improved on a GPU by increasing both $N$ and $M$. Due to the relevance of fully connected spin systems in real applications,  we have considered in this evaluation spin-glass Hamiltonians on complete graphs with $N = 2^n$ ($7 \leq n \leq  12$), and we performed simulated annealing once ($M = 1$) or $128$ times ($M = 128$), where each annealing trial consisted of executing $L = 20000$ Markov chain steps.
In the results summarized in Table \ref{tab:annealing}, each value represents the total annealing time per replica (i.e., the total annealing time $t_A$ divided by $M$) measured in milliseconds. Let us note that the algorithms implemented on a GPU have significantly shorter execution times compared to their CPU implementation counterparts.
Although the annealing time per trial for the Glauber dynamics, SCA, and $\vep$-SCA are approximately the same when implemented on a GPU, such an observation can be compensated by the fact that, as we will see in the following discussion, the $\vep$-SCA and the SCA outperform Glauber dynamics in terms of the success rate of finding a ground state.

 In order to understand the dependence of the performance of the algorithms with respect to certain variables, such as the system size $N$ and the number of Markov Chain steps $L$, we performed simulations on a GPU specifically for the problem of minimization of a spin-glass Hamiltonian on a complete graph considering a fixed number of replicas $M = 1024$ and a cooling schedule such as in equation (\ref{eq:cooling}) with initial temperature $T_\text{init} = 1000$ and final temperature $T_\text{fin} = 0.05$. In Table \ref{tab:comparison}, corresponding to each algorithm, we analyze its success rate $\rho_\text{min}$ of hitting the smallest energy reached among all algorithms, the average minimum energy $H_{\text{mean}}^{\text{(min)}}$ obtained, and total annealing time per replica $t_A / M$. As expected, we observe that for all algorithms, the success rates increase and the average minimum energy decrease as the number of Markov chain steps assume larger values. In that respect, the values of $\rho_\text{min}$ and $H_{\text{mean}}^{\text{(min)}}$  for SCA and $\vep$-SCA (with $\vep = 0.8$) are consistently better compared to the Glauber dynamics. Furthermore, the performances associated with the Glauber dynamics and SCA rapidly decrease as we consider problems with a larger number of vertices, while the $\vep$-SCA keeps the highest performance among the three algorithms.
}  

The greater effectiveness of the $\vep$-SCA in reaching lower energy configurations in several observations is very intriguing due to the lack of any rigorous mathematical justification (at the moment) for that. Therefore, it raises several questions to be answered that bring us a new direction to be explored for the development of efficient algorithms for obtaining ground states of Ising Hamiltonians. Moreover, in future investigations, hardware-related questions should be addressed, such as memory bandwidth problems and efficient implementation of parallel spin-flip algorithms on a GPU, similar to those explored in \cite{CM22}.

%\section*{Declarations}
%This work was supported by JST CREST Grant Number JP22180021, Japan. All data generated or analyzed during this study are %included in this published article. The authors have no relevant financial or non-financial interests to disclose.

%\section*{Data availability statement}

\section*{Acknowledgements}
This work was supported by JST CREST Grant Number JP22180021, Japan.
We are grateful to 
the following members for continual encouragement and stimulating discussions: 
Masato Motomura  from Tokyo Institute of Technology; 
Shinya Takamaeda-Yamazaki from the University of Tokyo; 
Hiroshi Teramoto from Kansai University; 
Takashi Takemoto, Takuya Okuyama, Normann Mertig and Kasho Yamamoto from 
Hitachi Ltd.; 
Masamitsu Aoki and Suguru Ishibashi from the Graduate School of 
Mathematics at Hokkaido University; Hisayoshi Toyokawa from Kitami Institute of Technology; 
Yuki Ueda from Hokkaido University of Education.

\bibliographystyle{plain}
\bibliography{bibliography}
\end{document}